\documentclass{article}

\usepackage{mathptmx}
\usepackage{helvet}
\usepackage{courier}
\usepackage{type1cm}
\usepackage{graphicx}        % standard LaTeX graphics tool
\usepackage[bottom]{footmisc}% places footnotes at page bottom
\usepackage{amsmath}
\usepackage{amsfonts}
\usepackage{mathrsfs}
\usepackage{dsfont}
\usepackage{amssymb}
\usepackage{enumerate}
\usepackage{graphicx}

\usepackage{color}
\usepackage[dvipsnames]{xcolor}
\usepackage{longtable}
\usepackage{verbatim}
\usepackage{booktabs}

\usepackage{algorithm}
\usepackage{algorithmicx}
\usepackage{algcompatible}
\usepackage{algpseudocode}

\newcommand{\BB}{{\MR} }
\newcommand{\CRA }{h }
\newcommand{\PRA }{q }

\newtheorem{theorem}{Theorem}[section]
\newtheorem{lemma}[theorem]{Lemma}

\newtheorem{corollary}{Corollary}[section]

\newtheorem{example}[theorem]{Example}

\numberwithin{equation}{section}

\newcommand{\bea}{\begin{eqnarray*}}
\newcommand{\eea}{\end{eqnarray*}}
\newcommand{\be}{\begin{eqnarray}}
\newcommand{\ee}{\end{eqnarray}}

\newcommand{\vsp} {\vspace{0.3cm}}

\def\raL2{\stackrel{L_2}{\ra}}

\def\e1{\mathrm{e}}

\def\NN{\mathds{N}}
\def\RR{\mathds{R}}

\def\esph1{\hspace*{1cm}}
\def\esph2{\hspace*{2cm}}
\def\Arg{\mathrm{Arg}}

\def\CR{\mathsf{CR}}

\def\MR{\mathsf{MR}}

\def\vol{\mathrm{vol}}

\def\ra{\rightarrow}

\def\ma{\alpha}

\def\me{\epsilon}

\def\mg{\gamma}

\def\mp{\partial}

\def\SB{{\mathscr B}}
\def\SC{{\mathcal C}}

\def\SX{{\mathscr X}}

\def\0b{\mathbf{0}}
\def\1b{\mathbf{1}}

\def\xb{\mathbf{x}}
\def\Xb{\mathbf{X}}

\newcommand{\carre} {\mbox{}~\hfill\rule{2mm}{2mm}}
\newcommand{\fin}
{\mbox{}~\hfill{\lower-0.3ex\hbox{$\triangleleft$}}}

\begin{document}

%\title[Asymptotically optimal quasi-uniform nested designs]{Quasi-uniform designs with optimal and near-optimal uniformity constant}
\title{Quasi-uniform designs with optimal and near-optimal uniformity constant}

%\author[L. Pronzato]{Luc Pronzato}
%
%\address{CNRS, Universit\'e C\^ote d'Azur, I3S, France}
%\curraddr{}
%\email{pronzato@i3s.unice.fr}
%\thanks{}
%
%%    author two information
%\author[A. Zhigljavsky]{Anatoly Zhigljavsky}
%\address{Cardiff University, UK}
%\curraddr{}
%\email{ZhigljavskyAA@cardiff.ac.uk}
%\thanks{}

\author{Luc Pronzato\footnotemark[1]\ \ and Anatoly Zhigljavsky\footnotemark[2] }
\renewcommand{\thefootnote}{\fnsymbol{footnote}}
\footnotetext[1]{CNRS, Universit\'e C\^ote d'Azur, I3S, France, {\tt pronzato@i3s.unice.fr}}
\footnotetext[2]{Cardiff University, UK, {\tt ZhigljavskyAA@cardiff.ac.uk}}

%\footnotetext[3]{pronzato@i3s.unice.fr}

\renewcommand{\thefootnote}{\arabic{footnote}}

\maketitle

%    \subjclass is required.
%\subjclass[2020]{Primary 65D17, 05B30; secondary 65D15}

%\date{\today}

%\dedicatory{}

%    Abstract is required.
\begin{abstract} A design is a collection of distinct points in a given set $\SX$, which is assumed to be a compact subset of $\RR^d$, and
the mesh-ratio of a design is the ratio of its fill distance to its separation radius.
The uniformity constant of a sequence of nested designs is the smallest upper bound for the mesh-ratios of the designs. We derive
a lower bound on this uniformity constant and show that a simple greedy construction achieves this lower bound. We then extend
this scheme to allow more flexibility in the design construction.

\medskip
\noindent AMS subject classifications:  {Primary 65D17, 05B30; secondary 65D15}
\end{abstract}

\section{Introduction}
Let $\SX$ be a compact subset of $\RR^d$, for some $d\geq 1$, with ${\rm vol}(\SX)>0$.

Let $\| \cdot \| $ denote a norm, not necessarily the Euclidean  norm $\|\cdot\|_2$, on $\RR^d$.
The ball of radius $r$ and center $\xb$ is $\SB(\xb,r) = \{\xb' \in \RR^d: \| \xb'-\xb \| \leq r\}$. The volume of
the unit ball $\SB(\0b,1)$ is denoted by $V_d$. If  the norm $\| \cdot \| $ is  Euclidean, then $V_d=\pi^{d/2}/\Gamma(d/2+1)$.

A collection $\Xb_n=\{\xb_1,\ldots,\xb_n\}$ of $n$ distinct points in $\SX$ will be called an $n$-point design (in the modern literature on approximation theory,   designs are often called ``data sets'', see e.g. \cite{SchabackW2006, Wendland2005}).
We start with several definitions of well-known  characteristics of designs.

 FD, {\it the fill distance} (also known as mesh norm, covering radius, dispersion, or minimax-distance criterion), of the $n$-point design $\Xb_n$ for $\SX$  is
\bea
\CRA  (\Xb_n)=\CRA_\SX  (\Xb_n) := \sup_{\xb\in\SX}\min_{\xb_i\in\Xb_n} \|\xb-\xb_i\|\,, \ n\geq 1\,.
\eea
A design $\Xb_{n,FD  }^*$ will be called FD-optimal  if
$
\CRA_n^*:=\CRA  (\Xb_{n,FD  }^*)=\min_{\Xb_n\in\SX} \CRA  (\Xb_n) \,.
$
SR, {\it the separation radius} (also called packing radius or maximin-distance criterion), of $\Xb_n$  is
\bea
\PRA (\Xb_n) = \frac12\,\min_{\xb_i\neq \xb_j \in\Xb_n} \|\xb_i-\xb_j\| \,, \ n\geq 2\,.
\eea
A design $\Xb_{n,\PRA }^*$ will be called SR-optimal if
$
\PRA _n^*:=\PRA (\Xb_{n,SR }^*)=\max_{\Xb_n\in\SX} \PRA (\Xb_n) \,.
$
{\it The mesh-ratio} of $\Xb_n$ for $\SX$ is
\bea
\MR(\Xb_n)=\MR_\SX(\Xb_n) := \frac{\CRA_\SX  (\Xb_n)}{\PRA (\Xb_n)} \,, \ n\geq 2\,.
\eea

The mesh-ratio  provides a measure of how
uniformly points in $\Xb_n$ are distributed in $\SX$, see e.g. \cite[p.~573]{SchabackW2006} and \cite[p.~129]{fasshauer2007meshfree}; it  is sometimes called the uniformity constant of $\Xb_n$, see \cite{de2010stability}.
The mesh-ratio is commonly used to asses the stability of approximations constructed on the base of observations at $\xb_i \in \Xb_n$, see e.g.
\cite{SchabackW2006} and \cite[Chapter~12]{Wendland2005}.
According to Guideline 7.10 in \cite[p.~579]{SchabackW2006}, the best approximation error with the most stable system is achieved by using quasi-uniform designs (data sets) with the smallest mesh-ratio. The mesh-ratio is fundamental in estimation of stability of approximations through the approach involving the Lebesgue constant, see  \cite[Th.~1]{de2010stability} and \cite[Sect.~8.5]{Iske2018}. Moreover, the mesh-ratio plays an an important role in the derivation of upper-bounds on the quality of kernel approximations in the so-called `escape theorems', when the approximated function is less smooth than the kernel, see \cite{NarcowichWW2004, NarcowichWW2006} as well as \cite[Th.~1, p.~129]{fasshauer2007meshfree} and \cite[Th.~7.8]{SchabackW2006}.

%A design $\Xb_{n,\MR}^*$, $n\geq 2$, will be called $\MR$-optimal if
%\bea
%\MR_n^*:=\MR(\Xb_{n,\MR}^*)=\min_{\Xb_n\in\SX} \MR(\Xb_n) \,.
%\eea
%Note that $n$-point $\CRA  $- and $\PRA $- and $\MR$-optimal designs always exist as $\SX$ is compact; however these designs are not necessarily unique.

%{\Luc Do we really need the distinction between $\Xb_\infty$ and $\{ \Xb_n\}_{n=1}^\infty$? The statement ``In particular,  $\BB(\Xb_\infty) \geq 2\,$ for any $\Xb_\infty \subset \SX$." in Theorem~\ref{th:1} is maybe a bit awkward as it seems to repeat what is written just above.}
%\AZ{Niederrieter always talks about sequences of points. I wanted to make a link with them and kind of demonstrate that we are better.}

Let $\Xb_\infty=\{\xb_1,\xb_2,\ldots\} \subset \SX$ be a sequence of points in $\SX$. There is a one-to one correspondence between such point sequence $\Xb_\infty$  and the sequence  $\{ \Xb_n\}_{n=1}^\infty$ of nested designs $\Xb_n=\{\xb_1,\ldots,\xb_n\}$.
 A sequence $\{ \Xb_n\}_{n=1}^\infty$ of nested designs $\Xb_n$ in a compact  set $\SX\subset\RR^d$ is called quasi-uniform if there exists a constant $b<\infty$ such that $\MR(\Xb_n)\leq b$ for all $n$. The smallest such $b=\BB(\Xb_\infty)$ is called the \emph{uniformity constant} of the corresponding sequence of nested designs $\{ \Xb_n\}_{n=1}^\infty$.
 Quasi-uniform sequences of designs with small uniformity constants are the main sources of designs (point sets) in the meshless  (or ``mesh-free'') methods of computational mathematics; see e.g. \cite{fasshauer2007meshfree,SchabackW2006,Wendland2005}

 A sequence $\Xb_\infty^*=\{\xb_1^*,\xb_2^*,\ldots\}$ will be called $\MR$-optimal if its uniformity constant is minimal:
 \be \label{eq:MR_optimal}
 \BB(\Xb_\infty^*)= \min_{\Xb_\infty \subset \SX} \BB(\Xb_\infty)\, .
 \ee

It is well known that when $\SX$ is connected, $\MR(\Xb_n)\geq 1$ for any $n$-point design $\Xb_n$ in $\SX$ (as the $n$-balls $\SB(\xb_i,\CR(\Xb_n))$ must cover $\SX$).
One of the main results of the paper is Theorem~\ref{th:1} below, which states that in fact $\limsup_{n\to\infty}\, \MR(\Xb_n)\geq 2$ for any compact $\SX$ with positive volume. The proof is rather elementary but the result does not seem to be known. It implies in particular that the classical greedy packing algorithm is $\MR$-optimal.

 \begin{theorem} \label{th:1}
 For any sequence of nested designs $\Xb_n$ in a compact  set $\SX\subset\RR^d$ with ${\rm vol}(\SX)>0$, we have
\bea
\limsup_{n\to\infty}\, \MR(\Xb_n)\geq 2 \,.
\eea
In particular,  $\BB(\Xb_\infty) \geq 2\,$ for any $\Xb_\infty \subset \SX$.
%
% For any compact   $\SX \subset \RR^d$ with ${\rm vol}(\SX)>0$, any norm  $\| \cdot \| $ on  $\RR^d$ and any point sequence
%$\Xb_\infty \subset \SX$, we  have $$\limsup_n \MR(\Xb_n) \geq 2\,,$$ where  $\Xb_n$ are the designs consisting of the first $n$ points of $\Xb_\infty$.
\end{theorem}

Theorem~\ref{th:1} is proved in Section~\ref{sec:proof}.
The greedy-packing (or coffee-house) algorithm is presented in Section~\ref{S:Algo-greedy-packing}; it constructs  a sequence $\Xb_\infty $ with $\MR(\Xb_n)\leq 2$ for all $n \geq 2$ and hence $\BB(\Xb_\infty) = 2$. In Section~\ref{S:relaxed-greedy-packing}, we generalize the greedy-packing algorithm to the construction of other quasi-uniform sequences with bounded $\BB(\Xb_\infty)$. In Section~\ref{sec:finite_cand} we use the results of Section~\ref{S:relaxed-greedy-packing} to establish properties of an implementable version of the greedy-packing algorithm where, at every iteration, the next design point $\xb_{n+1}$ is chosen among a finite set of candidates $\SX_N\subset\SX$ rather than within the whole $\SX$. In Section~\ref{sec:boundary-phobic-GP} we consider a boundary-phobic version of greedy packing, which provides designs with worse (larger) mesh-ratio but better (smaller) fill distance.
%Algorithm~\ref{algo:greedy-packing}

\section{Proof of Theorem~\ref{th:1}} \label{sec:proof}
Before providing a proof of  Theorem~\ref{th:1}, we prove two simple  lemmas, both  of them  presenting independent interest.

%
%\begin{lemma}\label{L:connectedX} Any design $\Xb_n$ in a connected set $\SX\subset\RR^d$, $n\geq 2$, satisfies $\CRA  (\Xb_n)\geq \PRA (\Xb_n)$; that is, $\MR(\Xb_n)\geq 1$.
%\end{lemma}
%%
%\begin{proof}
%Suppose that $\CRA  (\Xb_n)=h<\PRA (\Xb_n)$. Then, the $n$ balls $\SB(\xb_i,h)$, $\xb_i\in\Xb_n$, do not intersect (as $h<\PRA (\Xb_n)$) and cover $\SX$ (as $h=\CRA  (\Xb_n)$). This is impossible since $\SX$ is connected.
%\carre
%\end{proof}

%-----------------------------------------------------------
\begin{lemma}\label{L:CR&PR-anydesign} For any design $\Xb_n$ in a compact  set $\SX\subset \RR^d$, we have
\bea
\left[\vol(\SX)/V_d\right]^{1/d}\, n^{-1/d} \leq \CRA  (\Xb_n)\,, \ n\geq 1\,.
\eea
Moreover, for any $m$ such that $n \geq  m \geq 2$, we have
\bea
\PRA (\Xb_n) \leq \left[\vol(\SX_0)/V_d\right]^{1/d}\, n^{-1/d}\,,
\eea
where $\SX_0=\SX\oplus \SB(\0b,\PRA (\Xb_{m}))$, $\Xb_{m} $ is a sub-design of $ \Xb_n$ consisting of $m$ points and  $\oplus$ denotes the Minkowski sum.
\end{lemma}
\noindent {\em Proof.} %\begin{proof}
The $n$ balls $\SB(\xb_i,\CRA  (\Xb_n))$ cover $\SX$; this yields the first inequality. The second inequality follows from $\PRA (\Xb_n)\leq \PRA (\Xb_{m})$, which implies that all the balls $\SB(\xb_i,\PRA (\Xb_n))$ are fully inside  $\SX_0$ ($i=1, \ldots, n$).
\carre \\ \vsp

%%-----------------------------------------------------------
%\begin{lemma}\label{L:CR&PR} Any quasi-uniform sequence of nested designs $\Xb_n$ with uniformity constant\index{uniformity constant} $b$ in a compact  set $\SX\subset\RR^d$ satisfies, for any $n_0\geq 2$,
%\bea
%\frac{\vol^{1/d}(\SX)}{V_d^{1/d}}\, n^{-1/d} \leq \CRA  (\Xb_n) \leq b\, \PRA (\Xb_n) \leq \frac{b\, \vol^{1/d}(\SX_0)}{V_d^{1/d}}\, n^{-1/d}\,, \ \forall n\geq n_0\,,
%\eea
%where $\SX_0=\SX\oplus \SB(\0b,\PRA (\Xb_{n_0}))$, with $\oplus$ denoting the Minkowski sum and $V_d$ the volume of the $d$-dimensional unit ball $\SB(\0b_d,1)$.
%\end{lemma}
%%
%\begin{proof}
%The $n$ balls $\SB(\xb_i,\CRA  (\Xb_n))$ cover $\SX$; this yields the leftmost inequality. $\PRA (\Xb_n)\leq \PRA (\Xb_{n_0})$ for all $n\geq n_0$, so that all the balls $\SB(\xb_i,\PRA (\Xb_n))$ are included in $\SX_0$. This yields the rightmost inequality. The middle one is simply the definition of a quasi-uniform sequence.
%\carre
%\end{proof}

Lemma~\ref{L:CR&PR-anydesign} has the following  consequence concerning  the rate of decrease of the fill distance  and separation radius of quasi-uniform sequences of nested designs.

\begin{corollary}\label{Coro:CR&PR-sequence}
For any quasi-uniform sequence of nested designs $\Xb_n$
with uniformity constant $\rho$
in a compact  set \mbox{$\SX\subset\RR^d$}, we have
\be
\label{eq:speed}
c_1\,n^{-1/d} \leq \CRA  (\Xb_n) \leq \rho\,\PRA (\Xb_n) \leq c_2\,n^{-1/d} \,, \ \forall n\geq 2\,,
\ee
where $c_1$ and $c_2$ are  some positive constants.
\end{corollary}

In the case of Euclidean norm, the statement of Corollary~\ref{Coro:CR&PR-sequence} is proved in \cite{Wendland2005}; see Proposition 14.1 and the discussion just after it.

%-----------------------------------------------------------
\begin{lemma}\label{L:PR-CR-n-n+1}
Let, for any given $n\in\NN$,   $\Xb_n=\{\xb_1,\ldots,\xb_n\} $ and  $\Xb_{n+1}' =\{\xb_1',\ldots,\xb_{n+1}'\}$ be arbitrary $n$-point  and $(n\!+\!1)$-point designs in $\SX$. Then
\bea
\PRA (\Xb_{n+1}') \leq \CRA  (\Xb_n) \,.
\eea
\end{lemma}
\noindent {\em Proof.} %\begin{proof}
Since the $n$ balls $\SB(\xb_i,\CRA  (\Xb_n))$ cover $\SX$, the pigeon-hole principle implies that at least one of them must contain at least two points $\xb_i'$ and $\xb_j'$ from $\Xb_{n+1}'$. Therefore, $\|\xb_i'-\xb_j'\|\leq 2\,\CRA  (\Xb_n)$, implying  $\PRA (\Xb_{n+1}') \leq \CRA  (\Xb_n)$.
\carre \\ \vsp

%-----------------------------------------------------------

%%-----------------------------------------------------------
%
%Note that Lemma~\ref{L:CR&PR-anydesign} applies to $\CRA  $- and $SD$-optimal designs.

%-----------------------------------------------------------
%\begin{lemma}\label{L:MR=2isoptimal}
%For any sequence of nested designs $\Xb_n$ in a compact  set $\SX\subset\RR^d$ satisfies
%\bea
%\limsup_{n\to\infty} \MR(\Xb_n)\geq 2 \,.
%\eea
%\end{lemma}

\medskip
\noindent {\it Proof of Theorem~\ref{th:1}}. Assume that $\limsup_{n\to\infty} \MR(\Xb_n) < 2$. This would yield that there exists $r<2$ and $n_0$ such that $ \MR(\Xb_n) \leq r$ for all $n \geq n_0$.

Consider all such $n \geq n_0$.
The definition of $\CRA  (\Xb_n)$ and $\MR(\Xb_n)$ imply the existence of $\xb_j\in\Xb_n$ such that $$\|\xb_{n+1}-\xb_j\|\leq\CRA  (\Xb_n) \leq r\,\PRA (\Xb_n)\, .$$ Therefore,
\bea
\PRA (\Xb_{n+1}) \leq (1/2)\,\min_{\xb_i\in\Xb_n}  \|\xb_{n+1}-\xb_i\| \leq (r/2)\, \PRA (\Xb_n)\,.
\eea
 This implies the exponential decrease of $\PRA (\Xb_n)$ to zero (as $n \to \infty$), which contradicts \eqref{eq:speed}.
\carre
%\medskip
%\vsp

\section{Construction of sequences of quasi-uniform designs }\label{S:Algo-greedy-packing0}

\subsection{Greedy packing}\label{S:Algo-greedy-packing}

Let us first describe the greedy-packing algorithm (called ``geometric greedy method'' in \cite{de2005near}), which
achieves the lower bound of Theorem~\ref{th:1} and hence constructs an $\MR$-optimal sequence of points $\Xb_\infty$ and nested designs $\{\Xb_n\}_{n=1}^\infty$.
This algorithm is sometimes called the ``coffee-house" algorithm, due to the analogy with the behavior of customers in large coffee shops, where new clients tend to seat as far as possible  from occupied tables \cite{Muller2007}.

\begin{algorithm}[ht]
{\small
\caption{(Greedy packing)}\label{algo:greedy-packing}
\begin{algorithmic}[1]
\Require $\SX$ compact subset of $\RR^d$,  $\xb_1\in\SX$.
\State set $n=1$, $\Xb_1=\{\xb_1\}$;
\State for {$n=1,2, \ldots$} do the following:
\State find $\xb_{n+1}\in\Arg\max_{\xb\in\SX} \min_{\xb_i\in\Xb_n} \|\xb-\xb_i\|$,
 \State set $\Xb_{n+1}=\Xb_n\cup\{\xb_{n+1}\}$.

\end{algorithmic}
}
\end{algorithm}

For arbitrary   $\xb_1\in\SX$ and any choice of $\xb_{n+1}\in\Arg\max_{\xb\in\SX} \min_{\xb_i\in\Xb_n} \|\xb-\xb_i\|$ at step~3, the sequence of designs $\Xb_n$ constructed by Algorithm~\ref{algo:greedy-packing} satisfies the following  property.
%; {\Luc see \cite[Lemma~5.1]{de2005near}}.
%\AZ{ \cite[Lemma~5.1]{de2005near} give inequality only, not the equality}

\begin{lemma} \label{th:2a} For all $n\geq 2$, the designs $\Xb_n$ generated by Algorithm~\ref{algo:greedy-packing} satisfy
$\PRA(\Xb_n)=\CRA(\Xb_{n-1}){ /2}$.
\end{lemma}
\noindent {\em Proof.} %\begin{proof}
The inequality $\PRA(\Xb_n)\geq \CRA(\Xb_{n-1}){ /2}$ is proved in \cite[Lemma~5.1]{de2005near} by induction on $n$; the equality is obtained by the same arguments. By the definition of $\xb_2$, we have $\PRA (\Xb_2)=\CRA  (\Xb_1)/2$.

Assume that $\PRA(\Xb_n)=\CRA(\Xb_{n-1})/2$ and consider $\PRA(\Xb_{n+1})$:
\bea
\PRA(\Xb_{n+1}) &=& \min\left\{\PRA(\Xb_n), (1/2)\,\min_{\xb_i\in\Xb_n}\|\xb_{n+1}-\xb_i\|\right\}  \\
&=& \min\left\{\PRA(\Xb_n), \CRA(\Xb_n)/2\right\} \\
&=& \min\left\{\CRA(\Xb_{n-1})/2,\CRA(\Xb_n)/2\right\} = \CRA(\Xb_n)/2 \,. \hspace{4cm} \carre
\eea
%\\ \vsp

\begin{theorem} \label{th:2}  For all $n \geq 2$, the designs $\Xb_n$ generated by Algorithm~\ref{algo:greedy-packing} satisfy
\bea
\CRA  (\Xb_n) \leq 2\,\CRA  _n^*\,, \;\;
\PRA (\Xb_n) \geq \frac12\, \PRA _n^*\, ,\;\;
\MR(\Xb_n) \leq 2 \,.
\eea
\end{theorem}

\noindent {\em Proof.} %\begin{proof}
By Lemma~\ref{L:PR-CR-n-n+1} applied to the designs  $\Xb_{n+1}$  and $\Xb_{n,FD  }^*$, we obtain $\PRA (\Xb_{n+1})\leq \CRA  _n^*$. Using Lemma~\ref{th:2a}, this gives  $\CRA  (\Xb_n)\leq 2\, \CRA  _n^*$. From Lemma~\ref{L:PR-CR-n-n+1} applied to the designs $\Xb_{n+1,SR }^*$ and $\Xb_{n}$ and Lemma~\ref{th:2a}, we obtain $\PRA _{n+1}^*\leq \CRA  (\Xb_n) = 2\, \PRA (\Xb_{n+1})$. Finally, $\MR(\Xb_{n+1})=\CRA  (\Xb_{n+1})/\PRA (\Xb_{n+1})\leq \CRA  (\Xb_n)/\PRA (\Xb_{n+1}) =2\, .$
\carre \\ \vsp

Theorem~\ref{th:2} may be deduced from  Theorem 2.2
in  \cite{Gonzalez85}, where  Algorithm~\ref{algo:greedy-packing} is used to minimize the maximum intercluster distance; see also \cite[Theorem 4.3]{har2011geometric}.
%{\Luc I think that our proof is much too similar to \cite[Lemma~5.1]{de2005near} to write the next sentences. But I agree with your suggestion too.} \AZ{Our  proof of Theorem 2 is similar but not really of Theorem 1. They missed equality.}  {\em Theorem~\ref{th:2} also follows from Lemma 5.1 in \cite{de2005near} and from Theorem~\ref{th:3} below, when $a=1$. We, however, provided a proof of Theorem~\ref{th:2} as we believe our proof is the simplest and uncovers the property of Algorithm~\ref{algo:greedy-packing}, formulated as Lemma~\ref{th:2a}, which has not been noticed in the publications cited above.} {\Luc I rather suggest the following:}
Theorem~\ref{th:2} also follows from Theorem~\ref{th:3} below.
However, we think that the proof provided above is interesting in itself, as the important role of Lemma~\ref{th:2a} uncovers the key property of Algorithm~\ref{algo:greedy-packing}.

Note that in Theorem~\ref{th:2} the choice of the norm in $\SX$ is irrelevant.
Moreover, $\SX$ does not have to be a subset of $\RR^d$; in particular, $\SX$ can be a discrete set as in the clustering problems considered in \cite{Gonzalez85}.

While the calculation of $\PRA(\Xb_n)$ is straightforward, $\CRA(\Xb_n)$ is difficult to compute when $\SX$ is a continuous set. Methods of computational geometry can sometimes be used \cite{pronzato2017minimax}, but are restricted to  low-dimensional spaces. The substitution of a finite set $\SX_N$ for $\SX$, with the $N$ points of $\SX_N$ suitably well spread over $\SX$, is often used in practice; see Section~\ref{sec:finite_cand} for the analysis of this version of Algorithm~\ref{algo:greedy-packing}.

For $d=1$ and $\SX=[0,1]$, Algorithm~\ref{algo:greedy-packing} initialized at $x_1=1/2$ is equivalent to the celebrated van der Corput sequence in base 2 in terms of the behaviour of $\CRA(\Xb_n)$, $\PRA(\Xb_n)$  and $\MR(\Xb_n)$; see \cite[p.~25]{Niederreiter92}. The regular pattern of $\MR(\Xb_n)$ observed in dimension 1 extends to dimension 2 with $\SX=[0,1]^2$ when $\|\cdot\|=\|\cdot\|_2$ and the algorithm is initialized at the center $(1/2,1/2)$. This is illustrated on the left panel of Figure~\ref{F:fig1}: $\MR(\Xb_n)$ takes two values only, 2 and $\sqrt{2}$. The detailed behaviour of the algorithm is as follows.

\begin{theorem}\label{Th:CH-in-the-square} For any $n\geq 5$, define $m=m(n)=\lfloor \log_2(\sqrt{n/2-1/4}-1/2)\rfloor$. Then the packing and covering performance of Algorithm~\ref{algo:greedy-packing} with $\|\cdot\|=\|\cdot\|_2$, initialized at the center $(1/2,1/2)$ of $\SX=[0,1]^2$, is as follows:
\bea
\begin{array}{llll}
\PRA(\Xb_n)=\mg_m \sqrt{2}/4   \,, \!\!&\!\! \CRA(\Xb_n)=\mg_m/2 \,, \!\!&\!\! \MR(\Xb_n) = \sqrt{2} \,, \!\!&\!\! \mbox{for } n=n_m  \,, \\
\PRA(\Xb_n)=\mg_m/4    \,, \!\!&\!\! \CRA(\Xb_n)=\mg_m/2 \,, \!\!&\!\! \MR(\Xb_n) = 2\,, \!\!&\!\! \mbox{for } n=n_m+1,\ldots,k_m-1 \,, \\
\PRA(\Xb_n)=\mg_m/4 \,, \!\!&\!\! \CRA(\Xb_n)=\mg_m \sqrt{2}/4 \,, \!\!&\!\! \MR(\Xb_n) = \sqrt{2} \,, \!\!&\!\! \mbox{for } n=k_m  \,, \\
\PRA(\Xb_n)=\mg_m \sqrt{2}/8\,, \!\!&\!\! \CRA(\Xb_n)=\mg_m \sqrt{2}/4 \,, \!\!&\!\! \MR(\Xb_n) = 2\,, \!\!&\!\! \mbox{for } n=k_m+1,\ldots,n_{m+1}-1 \,,
\end{array}
\eea
where $\mg_m=2^{-m}$, $n_m=(2^m+1)^2+4^m$ and $k_m=(2^{m+1}+1)^2$.
\end{theorem}

For the sake of brevity, we only give a sketch of the full proof. It is based on the self-replicating pattern of the construction. The first five points in $\SX=[0,1]^2$ correspond to the corners and the center of the square. This gives the initialization for the beginning of the initial cycle, indexed by $m=0$, with $m$ denoting the cycle number. Define the initialization of cycle $m$ as the replication of the initial design of cycle $0$ into $4^m$ squares of side length $\mg_m=2^{-m}$, which form a regular partition of $[0,1]^2$. The initial design for cycle $m$ has thus $n_m=(2^m+1)^2+4^m=2^{2m+1}+2^{m+1}+1$ points: $(2^m+1)^2$ of them form a regular grid of width $\mg_m$ (i.e., a $(2^m+1)^2$ full factorial design); the other $4^m$ points are the centers of the small squares. When moving to the next cycle, the algorithm first (\textit{i}) adds the midpoints of the sides of all small squares (in arbitrary order), then (\textit{ii}) adds the $4^{m+1}$ centers of the smaller squares created at previous phase. The number of points added during phase (\textit{i}) equals $\ell_m=(2^{m+1}+1)^2-[(2^m+1)^2+4^m]=2^{m+1}(2^m+1)$.
For any $n\geq 5$, the associated cycle number $m=m(n)$ is the unique integer satisfying $n_m \leq n < n_{m+1}$. As $n_m=2(2^m+1/2)^2+1/2$, this gives $m(n)=\lfloor \log_2(\sqrt{n/2-1/4}-1/2)\rfloor$.

\begin{example}\label{Ex:Ex1} We take $\SX=[0,1]^2$, $\|\cdot\|=\|\cdot\|_2$ and $\xb_1=(1/2,1/2)$. Algorithm~\ref{algo:greedy-packing} progressively imbeds regular grids in $\SX$.
The left panel of Figure~\ref{F:fig1} shows the evolution of $\MR(\Xb_n)$ as a function of $n=2,\ldots,85$; the right panel shows $\Xb_n$ for $n=80=k_2-1$.

\begin{figure}[ht!]
\begin{center}
\includegraphics[width=.49\linewidth]{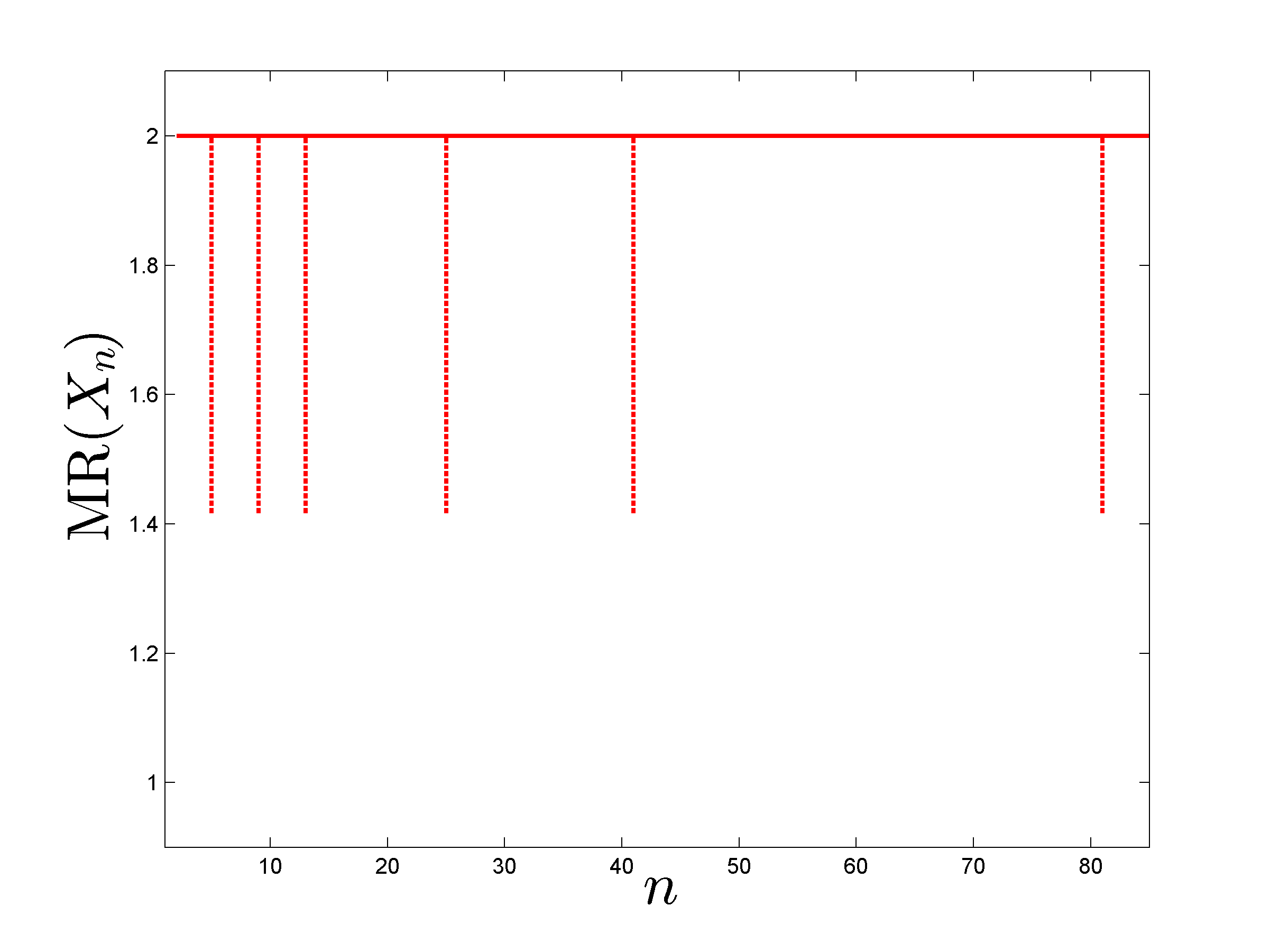}  \includegraphics[width=.49\linewidth]{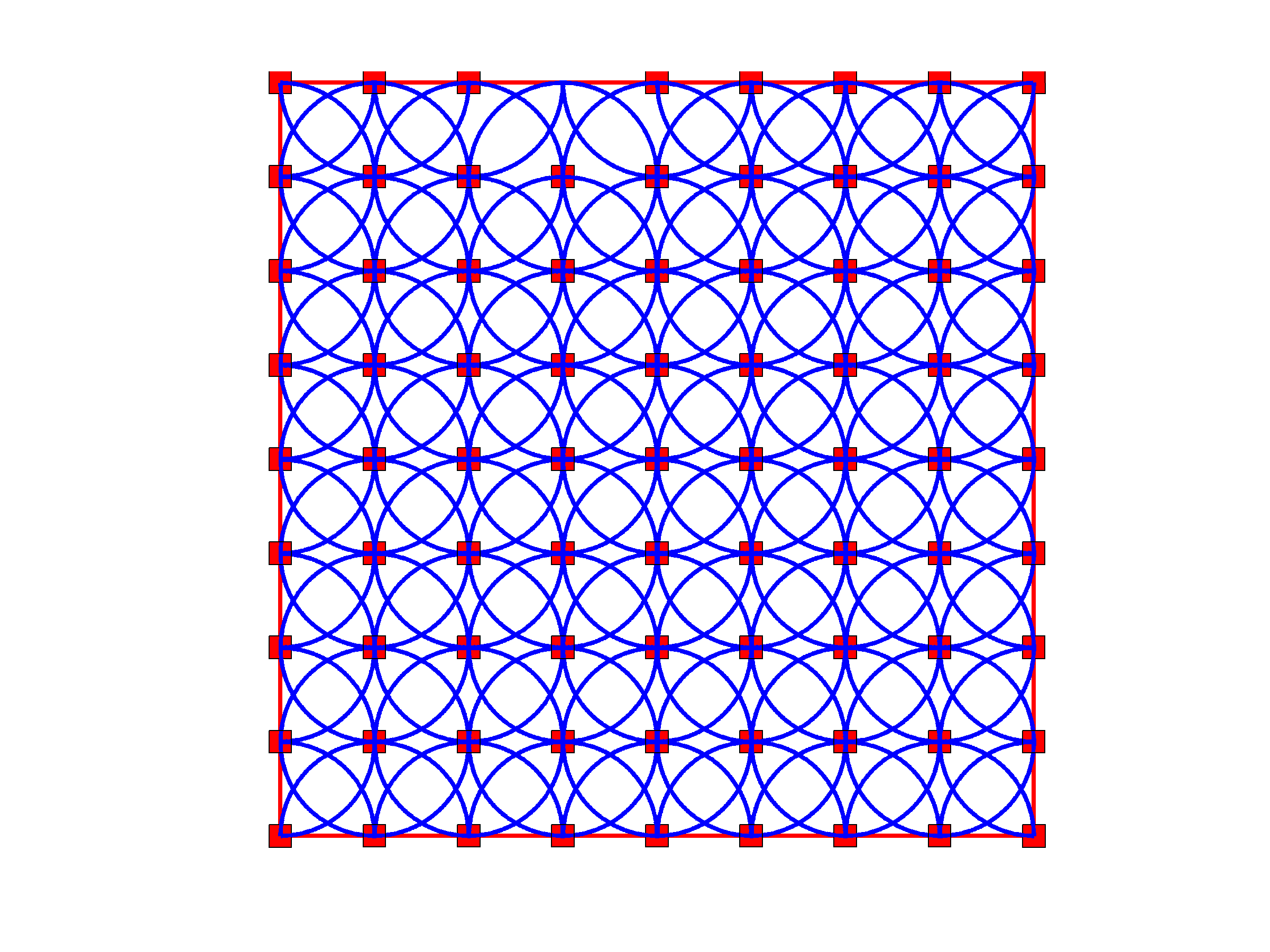}
\end{center}
\caption{\small Designs generated by Algorithm~\ref{algo:greedy-packing} in $\SX=[0,1]^2$ with $\|\cdot\|=\|\cdot\|_2$ and $\xb_1=(1/2,1/2)$. Left: $\MR(\Xb_n)$ for $n=2,\ldots,85$. Right: $\Xb_{80}$; the circles have radii $\CRA(\Xb_{80})= \mg_2/2 = 0.125$. }
\label{F:fig1}
\end{figure}

\end{example}

The regular pattern observed on $[0,1]^d$ for $d=1,2$ is maintained for $d=4$, and Algorithm~\ref{algo:greedy-packing} has the following behaviour in $[0,1]^4$.
\begin{theorem}\label{Th:CH-in-the-cube-d4} For any $n\geq 17$, define $m=m(n)$
as the unique integer satisfying $n_m \leq n < n_{m+1}$, with $n_m=(2^m+1)^4+2^{4m}$.
Then the packing and covering performance of Algorithm~\ref{algo:greedy-packing} with $\|\cdot\|=\|\cdot\|_2$, initialized at the center $(1/2,1/2,1/2,1/2)$ of $\SX=[0,1]^4$, is as follows:
\bea
\begin{array}{llll}
\PRA(\Xb_n)=\mg_m /2   \,, \!\!\!\!&\!\!\!\! \CRA(\Xb_n)=\mg_m \sqrt{2}/2 \,, \!\!\!\!&\!\!\!\! \MR(\Xb_n) = \sqrt{2} \,, \!\!\!\!&\!\!\!\! \mbox{ for } n=n_m  \,, \\
\PRA(\Xb_n)=\mg_m/(2\sqrt{2}) \,, \!\!\!\!&\!\!\!\! \CRA(\Xb_n)=\mg_m\sqrt{2}/2 \,, \!\!\!\!&\!\!\!\! \MR(\Xb_n) = 2\,, \!\!\!\!&\!\!\!\! \mbox{for } n=n_m+1,\ldots,n_m+\ell_m-1 \,, \\
\PRA(\Xb_n)=\mg_m/(2\sqrt{2}) \,, \!\!\!\!&\!\!\!\! \CRA(\Xb_n)=\mg_m /2 \,, \!\!\!\!&\!\!\!\! \MR(\Xb_n) = \sqrt{2} \,, \!\!\!\!&\!\!\!\! \mbox{ for } n=n_m+\ell_m  \,, \\
\PRA(\Xb_n)=\mg_m /4\,, \!\!\!\!&\!\!\!\! \CRA(\Xb_n)=\mg_m /2 \,, \!\!\!\!&\!\!\!\! \MR(\Xb_n) = 2\,, \!\!\!\!&\!\!\!\! \mbox{for } n=n_m+\ell_m+1,\ldots,n_{m+1}-1 \,,
\end{array}
\eea
where $\mg_m=2^{-m}$ and $\ell_m=6\times 2^{2m}\,(2^m+1)^2$.
\end{theorem}

The proof is omitted. Similarly to the 2-dimensional case treated in Theorem~\ref{Th:CH-in-the-square}, the construction follows a self-replicating pattern.
The first 17 points in $\SX=[0,1]^4$ are the 16 vertices and the center of $\SX$. This gives the initialization for the beginning of the cycle $m=0$,
which consists of the following two stages:
(\textit{i}) the algorithm chooses (in arbitrary order) all points with two coordinates equal to 1/2 and the other two coordinates in $\{0,1\}$; there are $2^2 \times \binom{4}{2}=24$ such points; (\textit{ii}) the algorithm chooses (in arbitrary order) points with one coordinate 1/2 and the other three in $\{0,1\}$ (there are $2^3 \times \binom{4}{1}=32$ such points), points with three coordinates 1/2 and one in $\{0,1\}$ (there are $2\times \binom{4}{3}=8$ such points) and points with coordinates in $\{1/4,3/4\}$ (there are 16 such points).

The initialization of cycle $m$ is defined as the replication of the initial design of cycle $0$ into $2^{4m}$ hypercubes of side length $\mg_m=2^{-m}$, which form a regular partition of $[0,1]^4$. The initial design for cycle $m$ has thus $n_m=(2^m+1)^4+2^{4m}$ points: $(2^m+1)^4$ of them form a regular grid of width $\mg_m$; the other $2^{4m}$ points are the centers of the small hypercubes. We thus have $2^{4m}$ replications of the initial 17-point initial design, but in smaller hypercubes. In each of them, the selections made by the algorithm are similar to those of the cycle $m=0$.

From the description above, we can observe that the design $\Xb_{n_m}$ re-scaled by a factor $2^{m}$ gives the integer lattice $\mathds{Z}_4$
truncated to $(i_1,i_2,i_3,i_4)\in\{0,\ldots,2^{m}\}^4$. Moreover, when $n=n_m+\ell_m$, the  design $\Xb_{n}$ re-scaled by  $2^{m+1}$ gives the
so-called checkerboard lattice $D_4$ (the subset of the integer lattice $\mathds{Z}_4$ consisting of quadruples whose sum is even), truncated to $(i_1,i_2,i_3,i_4)\in\{0,\ldots,2^{m+1}\}^4$; note that $D_4$ is the densest packing lattice in the 4-dimensional space \cite[p.~9]{ConwayS1999}.

%The regular behaviour of Algorithm~\ref{algo:greedy-packing} observed for $d=1,2$, and $d=4$ where the properly re-scaled design $\Xb_n$ oscillates between the integer point lattice and the checkerboard lattice, does not hold for other dimensions $d$.
The regular behaviour of Algorithm~\ref{algo:greedy-packing} observed for $d=1$, and $d=2$ and 4 where the properly re-scaled design $\Xb_n$ oscillates between the integer point lattice and the checkerboard lattice, does not hold for other dimensions $d$

%{\Luc \emph{The sentence suggests that the two types of lattices are present in dimension 1 and 2. I think it should be reformulated.}}
%However, as demonstrated in Section~\ref{sec:3d} for  $d=3$,  we may impose this  pattern   at the expense of increase  of $\MR(\Xb_n)$ for some $n$.

%---------------------------------------------------------------------------------
\subsection{Relaxed greedy packing}\label{S:relaxed-greedy-packing}

We consider now a generalization of Algorithm~\ref{algo:greedy-packing}, where the next point at a given iteration is not necessarily the furthest away from current design points, but is guaranteed to be far enough from them. The bounds obtained in Theorem~\ref{th:3} are worse than those in Theorem~\ref{th:2}; however, it can be shown that the relaxation introduced may improve the covering properties of the design sequence generated; see Section~\ref{sec:boundary-phobic-GP}.

\begin{algorithm}[ht!]
{\small
\caption{(Relaxed greedy packing)}\label{algo:relaxed-greedy-packing}
\begin{algorithmic}[1]
\Require $\SX$ compact subset of $\RR^d$,  $\xb_1\in\SX$, $a\in(0,1]$, $\ma_1,\ma_2,\ldots\in[a,1]$;
\State set $n=1$, $\Xb_1=\{\xb_1\}$;
\State for {$n=1,2, \ldots$} do the following:
\State  take any $\xb^\prime$ such that $\min_{\xb_i\in\Xb_n}\|\xb^\prime-\xb_i\|\geq \ma_n\,\CRA(\Xb_n)$ and set $\xb_{n+1}=\xb^\prime$;
 \State set $\Xb_{n+1}=\Xb_n\cup\{\xb_{n+1}\}$.
\end{algorithmic}
}
\end{algorithm}

At step~3, the choice of $\xb_{n+1}$ is arbitrary provided it satisfies the condition indicated. Due to this flexibility, several existing algorithms form particular cases of Algorithm~\ref{algo:relaxed-greedy-packing}, which is fact defines a whole family of algorithms. In particular, one may first select
$\xb^*\in\Arg\max_{\xb\in\SX}\min_{\xb_i\in\SX} \|\xb-\xb_i\|$ and then take any point $\xb_{n+1}\in\SB(\xb^*,(1-\ma_n)\CRA(\Xb_n))$.

\begin{theorem}\label{th:3} For all $n\geq 2$, the designs $\Xb_n$ generated by any version of Algorithm~\ref{algo:relaxed-greedy-packing} satisfy
\bea
\CRA(\Xb_n) \leq\frac{2}{a}\,\CRA_n^*\,, \;\;
\PRA(\Xb_n) \geq \frac{a}{2}\, \PRA_n^*\,, \;\;
\MR(\Xb_n) \leq \frac{2}{a}\,.
\eea
\end{theorem}

\noindent {\em Proof.} %\begin{proof}
We first prove by induction that for all $n\geq 2$, $\PRA(\Xb_n)\geq (a/2)\,\CRA(\Xb_{n-1})$.

For $n=2$, by construction we have $\PRA(\Xb_2) \geq (\ma_1/2)\,\CRA(\Xb_1)\geq (a/2)\, \CRA(\Xb_1)$.

Assume that $\PRA(\Xb_n)\geq (a/2)\,\CRA(\Xb_{n-1})$ and consider $\PRA(\Xb_{n+1})$.
The induction assumption gives
\bea
\PRA(\Xb_{n+1}) &=& \min\left\{\PRA(\Xb_n), (1/2)\,\min_{\xb_i\in\Xb_n}\|\xb_{n+1}-\xb_i\|\right\}  \\
&\geq& \min\left\{\PRA(\Xb_n),(\ma_n/2)\, \CRA(\Xb_n)\right\} \\
&\geq& \min\left\{(a/2)\,\CRA(\Xb_{n-1}),(a/2)\, \CRA(\Xb_n)\right\} = (a/2)\, \CRA(\Xb_n) \,.
\eea

The inequality proved by induction implies
$$
\MR(\Xb_n)=\CRA(\Xb_n)/\PRA(\Xb_n)\leq \CRA(\Xb_{n-1})/\PRA(\Xb_n)\leq 2/a\,.
$$
Next, by Lemma~\ref{L:PR-CR-n-n+1}, $\CRA(\Xb_{n-1})\geq \PRA_n^*$, and therefore
$\PRA(\Xb_n)\geq (a/2)\,\CRA(\Xb_{n-1}) \geq (a/2)\,\PRA_n^*$.
The same lemma implies $\CRA_{n-1}^*\geq \PRA(\Xb_n)\geq (a/2)\,\CRA(\Xb_{n-1})$.
\carre \\ \vsp

Theorem~\ref{th:2} follows from Theorem~\ref{th:3} by taking $a=1$. As in Theorem~\ref{th:2}, the choice of the norm in $\SX$ is irrelevant
and $\SX$ does not have to be a subset of $\RR^d$. In the next section, Theorem~\ref{th:3} is used for assessing properties of an easily implementable  version of Algorithm~\ref{algo:greedy-packing}, where $\xb_{n+1}$ at step 3 is chosen from a finite set.

\subsection{Greedy packing for a finite candidate set} \label{sec:finite_cand}

Consider a version of Algorithm~\ref{algo:greedy-packing} where $\xb_{n+1}$ is chosen among a finite set of candidates $\SX_N\subset\SX$ rather than from the whole $\SX$.  This assumption
makes the implementation of Algorithm~\ref{algo:greedy-packing} much simpler but  naturally deteriorates its performance.
Such implementation of Algorithm~\ref{algo:greedy-packing} can be considered as a special case of Algorithm~\ref{algo:relaxed-greedy-packing}, and hence, as we show below in Theorem~\ref{Th:PerfGreedypackingGrid}, its performance over entire $\SX$ can be assessed.
Note that the total number of iterations must be  smaller than $N$, the number of candidate points: indeed, for $n\geq N$, the algorithm degenerates  as several points necessarily coincide in $\Xb_{N+j}$, $j\geq 1$.

\begin{lemma}\label{L:discreteCR}
For any $n$-point design $\Xb_n$ and any $N$-point set $\SX_N\subset\SX$ we have
\bea
\CRA_{\SX_N}(\Xb_n) \leq \CRA_\SX(\Xb_n) \leq  \CRA_{\SX_N}(\Xb_n)+\CRA_\SX(\SX_N) \,.
\eea
\end{lemma}
\noindent {\em Proof.} %\begin{proof}
The inequality $\CRA_{\SX_N}(\Xb_n) \leq \CRA_\SX(\Xb_n)$ follows from $\SX_N\subset\SX$. Next, denoting $\SX_N=\{\xb^{(1)},\ldots,\xb^{(N)}\}$, we have
\bea
\CRA_\SX(\Xb_n) &\leq& \sup_{\xb\in\SX}\min_{\xb_i\in\Xb_n} \min_{\xb^{(j)}\in\SX_N} \left(\|\xb-\xb^{(j)}\| + \|\xb^{(j)}-\xb_i\|\right) \\
&=& \sup_{\xb\in\SX} \left[ \min_{\xb^{(j)}\in\SX_N} \left(\|\xb-\xb^{(j)}\| + \min_{\xb_i\in\Xb_n} \|\xb^{(j)}-\xb_i\|\right) \right] \\
&\leq& \sup_{\xb\in\SX} \left[ \min_{\xb^{(j)}\in\SX_N} \|\xb-\xb^{(j)}\| + \max_{\xb^{(j)}\in\SX_N} \min_{\xb_i\in\Xb_n} \|\xb^{(j)}-\xb_i\| \right] \\
&=& \sup_{\xb\in\SX} \left[ \min_{\xb^{(j)}\in\SX_N} \|\xb-\xb^{(j)}\| + \CRA_{\SX_N}(\Xb_n) \right] \\
&=& \CRA_\SX(\SX_N) + \CRA_{\SX_N}(\Xb_n) \,. \hspace{6.5cm} \carre
\eea
%\carre \\ \vsp

\begin{theorem}\label{Th:PerfGreedypackingGrid} When Algorithm~\ref{algo:greedy-packing} uses a finite set of candidates $\SX_N\subset\SX$ and $n < N$, its performance satisfies
\be\label{PerfGreedypackingGrid}
\begin{array}{lcll}
\CRA_\SX(\Xb_n) &\leq& (2/\ma_n)\,\CRA_n^*\,, \ &\forall n\geq 1\,,\\
\PRA(\Xb_n) &\geq& (\ma_n/2)\, \PRA_n^*\,, \ &\forall n\geq 2\,,\\
\MR_\SX(\Xb_n) &\leq& 2/\ma_n\,, &\forall n\geq 2\,,
\end{array}
\ee
with $\ma_n=1-\CRA_\SX(\SX_N)/\CRA_\SX(\Xb_n)$.
\end{theorem}

\noindent {\em Proof.} %\begin{proof}
Denote $\me=\CRA_\SX(\SX_N)$, so that Lemma~\ref{L:discreteCR} gives $\CRA_{\SX_N}(\Xb_n) \leq \CRA_\SX(\Xb_n) \leq  \CRA_{\SX_N}(\Xb_n)+\me$. At step~3 of Algorithm~\ref{algo:greedy-packing}, we have
\bea
\min_{\xb_i\in\Xb_n} \|\xb_{n+1}-\xb_i\|=\CRA_{\SX_N}(\Xb_n)\geq \CRA_\SX(\Xb_n) - \me = \ma_n\,\CRA_\SX(\Xb_n) \,,
\eea
with $\ma_n=1-\me/\CRA_\SX(\Xb_n)$. Since $\CRA_\SX(\Xb_n)$ is non-increasing with $n$, $\ma_n$ is non-increasing too (it reaches zero when $\Xb_n$ has exhausted $\SX_N$, that is, when $k=N$). Theorem~\ref{th:3} with $\ma_n$ substituted for $a$ implies \eqref{PerfGreedypackingGrid}.
\carre \\ \vsp

As we do not know $\CRA_\SX(\Xb_n)$ and thus $\ma_n$, we can use the inequality $\CRA_\SX(\Xb_n) \geq \CRA_{\SX_N}(\Xb_n)$, which gives $\ma_n\geq a_n=1-\CRA_\SX(\SX_N)/\CRA_{\SX_N}(\Xb_n)$. The inequalities \eqref{PerfGreedypackingGrid} then remain true with $a_n$ substituted for $\ma_n$, as long as $a_n>0$.

A result similar to Theorem~\ref{PerfGreedypackingGrid} holds when the performance of Algorithm~\ref{algo:greedy-packing} is evaluated on a finite set $\SX_{N'}\supset\SX_N$ instead of $\SX$: we simply substitute $\SX_{N'}$ for $\SX$ and $\ma_n=1-\CRA_{\SX_N'}(\SX_N)/\CRA_{\SX_N'}(\Xb_n)$ is evaluated easily.

\subsection{Boundary-phobic greedy packing} \label{sec:boundary-phobic-GP}

Versions of the greedy packing algorithm that enforce boundary avoidance have been proposed in \cite{NogalesPR2021,ShangA2020}. There, at iteration $n\geq 2$, the next point $\xb_{n+1}$ is chosen in $\Arg\max_{\xb\in\SX} D_\beta(\xb,\Xb_n,\SX)$,
where
\be\label{Dbeta}
D_\beta(\xb,\Xb_n,\SX) = \min\left\{\min_{\xb_i\in\Xb_n} \|\xb-\xb_i\|,\, \beta\,d(\xb,\mp\SX)\right\}\,, \ \beta\in(0,\infty)\,,
\ee
with $d(\xb,\mp\SX)$ the distance from $\SX$ to the boundary of $\SX$. Note that this quantity is easily determined if $\SX$ has a simple shape, like a hypercube or a ball, but may be difficult to evaluate otherwise. For $\beta=\infty$, we define $D_\infty(\xb,\Xb_n,\SX)=\min_{\xb_i\in\Xb_n} \|\xb-\xb_i\|$ by continuity; the algorithm then coincides with Algorithm~\ref{algo:greedy-packing}. For $\beta=1$, $\xb_{n+1}$ is the center of (one of) the largest ball included in $\SX$ and not intersecting $\Xb_n$. For $\beta=2$, the algorithm corresponds to a greedy method for the solution of the traditional packing problem, for which the $n$ balls do not intersect and are constrained to be fully inside $\SX$. For $\beta>2$, the larger $\beta$ is, the more the balls are allowed to overshoot $\SX$, with their centers remaining inside $\SX$. When $\SX=[0,1]^d$ and $\|\cdot\|=\|\cdot\|_2$, the value $\beta=2\,\sqrt{2d}$ is recommended in \cite{ShangA2020}, while \cite{NogalesPR2021} recommends to let $\beta$ depend on the targeted number $n_{\max}$ of design points and suggests taking
\bea
\beta= \beta(n_{\max},d)=\frac{d}{2\,(n_{\max} V_d)^{-1/d}}-\sqrt{d} \,,
\eea
with $V_d=\pi^{d/2}/\Gamma(d/2+1)$. Both references illustrate the interest of using $\beta<\infty$ instead of Algorithm~\ref{algo:greedy-packing} in terms of fill distance $\CRA(\Xb_n)$. As shown below, for $\SX=[0,1]^d$ the boundary-phobic version of greedy packing becomes a particular case of Algorithm~\ref{algo:relaxed-greedy-packing}.

\begin{theorem}\label{Th:boundary-phobic=Algorithm2}
For $\SX=[0,1]^d$ and $\|\cdot\|=\|\cdot\|_2$, the boundary-phobic algorithm that chooses $\xb_{n+1}$ in $\Arg\max_{\xb\in\SX} D_\beta(\xb,\Xb_n,\SX)$ at iteration $n$, with
$D_\beta(\xb,\Xb_n,\SX)$ defined by \eqref{Dbeta} and $\beta\in(0,\infty)$, forms a particular instance of Algorithm~\ref{algo:relaxed-greedy-packing} with $\ma_n=a=1/(1+\sqrt{d}/\beta)$.
\end{theorem}

\noindent {\em Proof.} %\begin{proof}
Let $\SX=[0,1]^d$ and $\beta\in(0,\infty)$, $r_n=D_\beta(\xb_{n+1},\Xb_n,\SX)=\max_{\xb\in\SX} D_\beta(\xb,\Xb_n,\SX)$. Any $\xb\in\SX$ satisfies at least one of the two inequalities
\bea
\min_{\xb_i\in\Xb_n}\|\xb-\xb_i\| \leq r_n \,, \ d(\xb,\mp\SX) \leq r_n/\beta \,.
\eea
This implies that $\SX\setminus \{\xb\in\RR^d: d(\xb,\mp\SX) \leq r_n/\beta\} \subset \cup_{i=1}^k \SB(\xb_i,r_n)$. The inequalities $r_n\leq \beta\, d(\xb_{n+1},\mp\SX)\leq \beta/2$ imply that $2\,r_n/\beta\leq 1$, and the set $\SX\setminus \{\xb\in\RR^d: d(\xb,\mp\SX) \leq r_n/\beta\}$ is a hypercube $\SC_n$ with side length $1-2\, r_n/\beta$. This hypercube is covered by the $n$ balls $\SB(\xb_i,r_n)$, implying that
\bea
\CRA_\SX(\Xb_n) &=& \sup_{\xb\in\SX} \min_{\xb_i\in\Xb_n} \|\xb-\xb_i\| \\
&\leq& \sup_{\xb\in\SX}\left[\inf_{\xb'\in\SC_n} \left(\|\xb-\xb'\|+\min_{\xb_i\in\Xb_n} \|\xb'-\xb_i\|\right)\right] \\
&\leq& \sup_{\xb\in\SX} \inf_{\xb'\in\SC_n} \|\xb-\xb'\| + r_n \leq \sqrt{d}\,(r_n/\beta) + r_n \,. \\
\eea
%$\CRA(\Xb_n)\leq \sqrt{d}\,(r_n/\beta) + r_n$.
Since, by definition, $r_n \leq \min_{\xb_i\in\Xb_n}\|\xb_{n+1}-\xb_i\|$, we have
\bea
\min_{\xb_i\in\Xb_n}\|\xb_{n+1}-\xb_i\| \geq \frac{\CRA_\SX(\Xb_n)}{1+\sqrt{d}/\beta} \,,
\eea
and the algorithm is a particular instance of Algorithm~\ref{algo:relaxed-greedy-packing} with $\ma_n=a=1/(1+\sqrt{d}/\beta)$.
\carre \\ \vsp

Theorem~\ref{Th:boundary-phobic=Algorithm2} implies that the performance of this algorithm satisfies the bounds indicated in Theorem~\ref{th:3}.

\begin{example} We take $\SX=[0,1]^2$, $\|\cdot\|=\|\cdot\|_2$ and $\beta=4$. The left panel of Figure~\ref{F:fig2} shows the evolution of $\MR(\Xb_n)$ as a function of $n=2,\ldots,80$ when $\Xb_n$ is generated by $\xb_{n+1}\in\Arg\max_{\xb\in\SX} D_\beta(\xb,\Xb_n,\SX)$ with $\xb_1=(1/2,1/2)$; the upper bound $2(1+\sqrt{d}/\beta)$ on $\MR(\Xb_n)$ is indicated by a horizontal line. The right panel presents $\Xb_{80}$: comparison with the right panel of Figure~\ref{F:fig1} shows that boundary avoidance has significantly reduced $\CRA(\Xb_n)$. This reduction is obtained at the detriment of $\MR(\Xb_n)$ for some $\Xb_n$, as illustrated by the left panels of the two figures (note, however, that $\MR(\Xb_{80})< 2$ on Figure~\ref{F:fig2}).

\begin{figure}[ht!]
\begin{center}
\includegraphics[width=.49\linewidth]{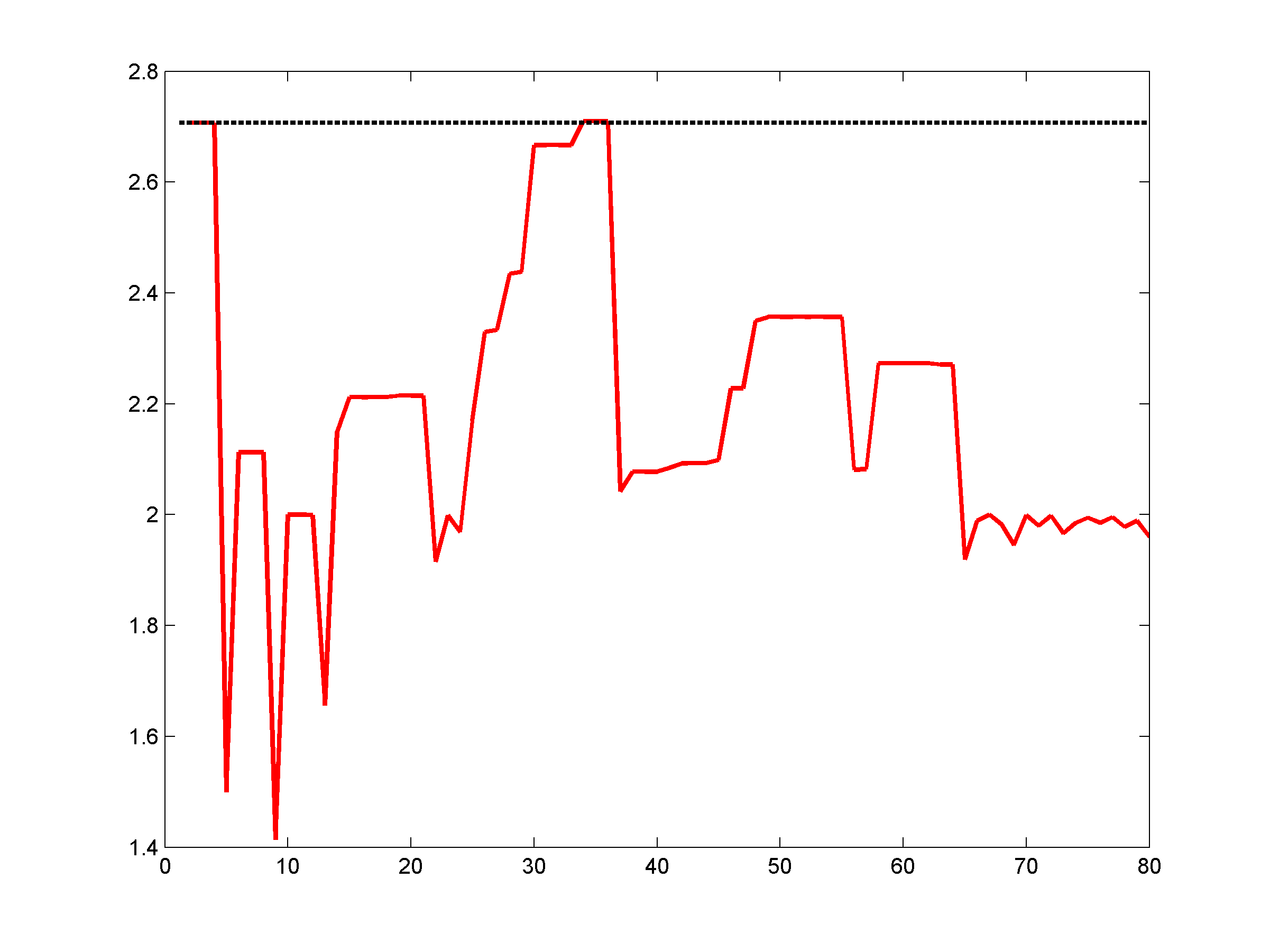} \includegraphics[width=.49\linewidth]{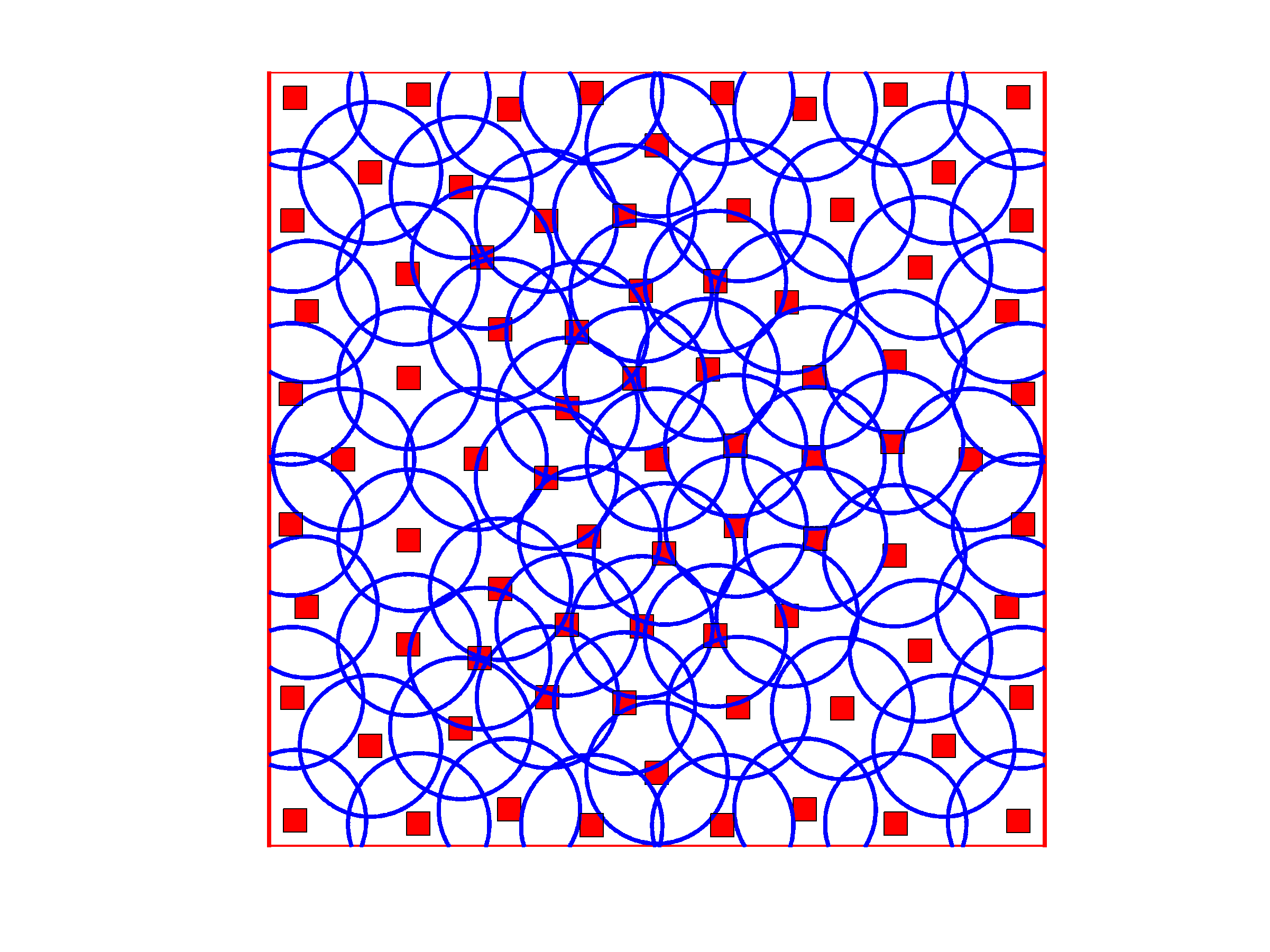}
\end{center}
\caption{\small Designs generated by $\xb_{n+1}\in\Arg\max_{\xb\in\SX} D_4(\xb,\Xb_n,\SX)$ in $\SX=[0,1]^2$ with $\|\cdot\|=\|\cdot\|_2$ and $\xb_1=(1/2,1/2)$. Left: $\MR(\Xb_n)$ for $n=2,\ldots,80$; the horizontal line indicates the upper bound $2(1+\sqrt{d}/\beta)$ on $\MR(\Xb_n)$. Right: $\Xb_{80}$; the circles have radii $\CRA(\Xb_n)=0.0913$.  }
\label{F:fig2}
\end{figure}

\end{example}

\bibliographystyle{amsplain}

%\bibliography{mesh_ratio,springer_luc}
%\bibliography{mesh_ratio}
\providecommand{\bysame}{\leavevmode\hbox to3em{\hrulefill}\thinspace}
\providecommand{\MR}{\relax\ifhmode\unskip\space\fi MR }
\providecommand{\MRhref}[2]{%
  \href{http://www.ams.org/mathscinet-getitem?mr=#1}{#2}
}
\providecommand{\href}[2]{#2}

\end{document}